\def\deffp<#1>{#1-bit floating-point}
\def\fp<#1>{{fp#1}}
\def\posit<#1,#2>{posit$\langle$#1,#2$\rangle$}
\def\p<#1>{posit#1}

\def\low<#1>{{posit#1}}
\def\working<#1>{{posit#1}}
\def\high<#1>{{posit#1}}

\documentclass{amsart}

\theoremstyle{definition}

\theoremstyle{remark}

\numberwithin{equation}{section}

\usepackage{tikz}
\usepackage{pgfplots}
\pgfplotsset{compat=1.15}
\usepackage{mathrsfs}
\usetikzlibrary{arrows}
\usepackage{amsmath, amsfonts}
\usepackage{xcolor}
\usepackage{hyperref,url} 
\hypersetup{
     colorlinks=true,
     linkcolor=black,
     filecolor=black,
     citecolor = black,      
     urlcolor=black, 
     }

\usepackage{xfrac}

\definecolor{sign}{RGB}{97, 160, 225}
\definecolor{fraction}{RGB}{255, 227, 116}   
\definecolor{fracrat}{RGB}{200, 170, 30}
\definecolor{expo}{RGB}{134,200,79}
\definecolor{regime}{RGB}{255,122,35}

\makeatletter
\renewcommand*\env@matrix[1][c]{\hskip -\arraycolsep
  \let\@ifnextchar\new@ifnextchar
  \array{*\c@MaxMatrixCols #1}}
\makeatother

\usepackage{algorithm}
\usepackage{algpseudocode}
\algrenewcommand\algorithmicrequire{\textbf{Input:}}
\algrenewcommand\algorithmicensure{\textbf{Output:}}

\usepackage{outlines}
\usepackage{enumitem}
\setenumerate[1]{label=\arabic*.}
\setenumerate[2]{label=(\alph*).}
\setenumerate[3]{label=\roman*.}

\def\num#1{\numx#1}\def\numx#1e#2{{#1}\mathrm{e}{#2}}

\begin{document}
\title{Iterative Refinement with Low-Precision Posits}
\author{James Quinlan}
\address{Department of Computer Science\\
                University of Southern Maine\\
                Portland, Maine}
\email{james.quinlan@maine.edu}
\thanks{Research partially supported by the Maine Economic Improvement Fund.}
\author{E. Theodore L. Omtzigt}
\address{Stillwater Supercomputing Inc.\\
		El Dorado Hills, California}
\email{tomtzigt@stillwater-sc.com}


\date{}
\dedicatory{}

\begin{abstract}
This research investigates using a mixed-precision iterative refinement method using posit numbers instead of the standard IEEE floating-point format. The method is applied to solve a general linear system represented by the equation $Ax = b$, where $A$ is a large sparse matrix. Various scaling techniques, such as row and column equilibration, map the matrix entries to higher-density regions of machine numbers before performing the $O(n^3)$ factorization operation. Low-precision LU factorization followed by forward/backward substitution provides an initial estimate.  The results demonstrate that a 16-bit posit configuration combined with equilibration produces accuracy comparable to IEEE half-precision (fp16), indicating a potential for achieving a balance between efficiency and accuracy.
\end{abstract}

\maketitle

\section{Introduction}\label{sed:intro}
Next-generation arithmetic is increasingly expected to address issues with power and performance in computing.
Both high-performance computing (HPC) and artificial intelligence (AI) fields are exploring nontraditional floating-point representations and arithmetic \cite{haidar2018design,kharya2020tensorfloat,lindquist2020improving,svyatkovskiy2017training}.  
Low precision formats such as IEEE's half-precision ({fp16}) and Google's brain float (\texttt{bfloat16}) can run at least twice the speed \cite{feldman2018fujitsu} and require only a fraction of the storage \cite{wu2021low}.
Such features, desirable for deep learning models \cite{gupta2015deep,svyatkovskiy2017training}, are also crucial in model simulation   \cite{palmer2014more}.
Manufacturers, including Intel \cite{intel2018bfloat16} and NVIDIA \cite{kharya2020tensorfloat}, are customizing hardware to support these alternate formats. 
Recently, hardware for the posit number system  \cite{gustafson2017beating} has also been introduced \cite{mallasen2022percival,suraksha2024silicon}.

{Our focus in this paper is to examine a low precision format in approximating the solution of a general linear system of equations} $A x = b$
where $A \in \mathbb{R}^{n \times n}$ is a nonsingular matrix and $b \in \mathbb{R}^n$.
Such systems are encountered in many scientific and industrial applications \cite{leon2020linear} and one of the most frequently occurring problems in computing \cite[p.16]{forsythe1967computer}.

In this work, we represent the entries of $A$ using a low-precision posit format to leverage the benefits of reduced power consumption and improved performance. 
However, we strive to maintain comparable levels of backward and forward errors, as accurate solutions are crucial in many applications spanning various domains, such as mathematical physics \cite{bailey2015high}, supercomputing \cite{he2001using}, optimization \cite{ma2015solving}, and molecular biology \cite{ma2017reliable}.

The paper is organized as follows: Section \ref{sec:ir} presents iterative refinement, mixed-precision iterative refinement, and some previous work on mixed-precision iterative refinement in the literature.  
Section \ref{sec:pns} covers the posit number system, gives an example, provides a table of properties, and compares it with other formats.
The details of our numerical experiments are presented in section \ref{sec:numx}, including the matrix test-suite, conversion algorithms to low-precision posits, and technical specifications for replication.  
Results are presented in Section \ref{sec:results}, and a summary is presented in Section \ref{sec:conclusions}.
We take $\| \cdot \| = \| \cdot \|_\infty$ to be the infinity norm for its practical error estimation \cite[p. 123-124]{golub1996matrix}.


\section{Iterative Refinement}\label{sec:ir}
Iterative refinement (Algorithm \ref{alg:ir}), introduced by Wilkinson  \cite{jwilkinson1963rounding}, consists of three steps to improve an approximate solution of a system of equations $Ax = b$ when $A$ is ill-conditioned. 
First, an initial solution $x_0$ is computed.
This step uses Gaussian elimination (i.e., $LU$), which requires $O(n^3)$ operations.
After an initial approximation, $x_{0}$, the method computes the residual $r_0 = b - Ax_{0}$ validates the initial solution's accuracy.  
A correction is then computed by  solving $Ac_{0} = r_{0}$ and the approximate is updated by $x_{1} = x_{0} + c_{0}$. 
The last two steps are iterated until convergence.  
For a detailed discussion of stopping criteria, refer to  \cite{barrett1994templates,demmel2006error}.
{However}, when $A$ is extremely ill-conditioned, iterative refinement may not lead to convergence, and additional modifications are needed.

\begin{algorithm}[ht!]
\caption{Iterative refinement}\label{alg:ir}
\begin{algorithmic}[1] 
\Require $A \in \mathbb{R}^{n \times n}$,  $b \in \mathbb{R}^n$, and maximum number of iterations $N$
\Ensure A vector  approximating solution $x$ to $Ax = b$ 
 \State Solve $Ax_0 = b$ 
 \For{$i = 0:N$ }
    \State $r_i = b - A x_i$ 
    \State Solve $A c_i = r_i$ 
    \State $x_{i+1} \gets x_i + c_i$
     \If{converged} 
        \State return $x_{i+1}$ 
    \EndIf
\EndFor
\end{algorithmic}
\end{algorithm}

Mixed-precision iterative refinement is a promising technique for achieving a high-precision solution to a linear equation by using low-precision arithmetic for complex computations. It aims to balance computational performance and solution accuracy. 
The main concept behind this technique is to carry out the $O(n^3)$ operations in Gaussian elimination using low-precision arithmetic and then storing the results in working precision. Subsequently, the approximation is refined iteratively by solving the correction equation using the high-precision residual, as detailed in the algorithm. 
In Step 6, we utilize the \textit{quire} memory register to defer rounding in our experiments. This iterative process enables the solution to attain the desired level of accuracy while minimizing the computational overhead associated with high-precision operations.

Several researchers have investigated mixed-precision iterative refinement. Langou et al. introduced a mixed-precision iterative refinement approach where the matrix factorization step is performed in single precision (32-bit). At the same time, the remaining operations are carried out in double precision (64-bit) \cite{baboulin2009accelerating,langou2006exploiting}.
Carson and Higham have recently explored versions of mixed-precision iterative refinement utilizing three different precisions \cite{carson2017new,carson2018accelerating}.
Haidar et al. have developed parallel implementations of these mixed-precision iterative refinement techniques \cite{haidar2017investigating}.
Higham et al. \cite{higham2019squeezing}  developed a mixed-precision GMRES-based iterative refinement algorithm using two-sided scaling to convert to half-precision format, which moves the elements of the largest magnitude close to the overflow threshold.
In \cite{buoncristiani2020evaluating}, researchers applied iterative refinement to symmetric matrices using various posit configurations (see Section \ref{sec:pns}).
Iterative refinement remains an active area of research, as evidenced by various contributions \cite{al2006lu,baboulin2009accelerating,carson2017new,carson2018accelerating,haidar2018design,higham1997iterative}.

\begin{algorithm}[h]
\caption{Mixed-precision iterative refinement}\label{alg:mixir}
\begin{algorithmic}[1] 
\Require $A \in \mathbb{R}^{n \times n}$ and  $b \in \mathbb{R}^n$
\Ensure A vector $\widehat{x}$ approximating solution $x$ to $Ax = b$ 
\State Factor $A = LU$ in low-precision 
\State Solve $LUx_0 = b$ in working-precision
\State $i \gets 0$
\State $r_0 = b - Ax_0$
\While{$||r_i|| > \epsilon \left(\|A\|_\infty \| x_i \|_\infty + \| b \|_\infty \right)$}
    \State $r_i = b - A x_i$ compute in high-preicsion
    \State Solve $A c_i = r_i$ 
    \State $x_{i+1} \gets x_i + c_i$
    \State $i \gets i+1$
\EndWhile
\end{algorithmic}
\end{algorithm}

%
\section{Posit Number System}\label{sec:pns}

The {\em Posit Number System} (PNS) is a set of machine numbers, denoted by $\mathbb{P}$, representing the real numbers $\mathbb{R}$ in binary developed by \cite{gustafson2017beating} in 2017.
Two parameters characterize a posit representation: (i) the total number of bits ${n}$ and (ii) the maximum number of exponent bits ${E}$.  
A \textit{standard} $n$-bit posit  ($E = 2$) \cite{posits2022}, is decoded as: 
\begin{align}\label{eqn:decode}
x &=  ((1 - 3s) + f) \times 2^{(1-2s)(4k + e + s)}
\end{align}
where the sign bit $s$ is $0$ if $x \ge 0$ and $1$ for $x < 0$. 
The exponent is $0 \le e \le E$. 
In the standard, the exponent $e \in \{0, 1, 2, 3\}$.  
For positive posits, (\ref{eqn:decode}) simplifies to, 
\[
x =  (1 + f) \times 2^{4k + e}
\]

The most notable distinction from IEEE-754 is the posit's dynamic \textit{regime}\footnote{There is always at least one regime bit, and for $n > 2$, there are at least two bits.  For $n>2$, $1 \le r \le n-1$, therefore $ -(n-1) \le k \le n-2$.} field.
While the exponent is a power-of-two scaling determined by the exponent bits, the regime is a power-of-$16$ scaling determined by (\ref{eqn:regime}) where $r$ is the number of identical bits in the run before an opposite terminating bit, ${r^c}$. 

\begin{align}\label{eqn:regime}
k &= \begin{cases}
-r & \quad \mbox{ if run of }\;  0 \mbox{'s    length } r \\ 
 r - 1 & \quad \mbox{ if run of }\; 1 \mbox{'s   length } r,
 \end{cases}
\end{align}
and the normalized fraction $f \in [0,1)$ given by
$$f = 2^{-m}\sum_{i=0}^{m-1} f_i 2^i$$ 
with $m$ the number of fraction bits ranging between $0$ and $\max(0,n-5)$ and $f_i \in \{0,1\}$ is the $i$th bit in fraction bitstring. 
The standard $n$-bit posit configuration is 
denoted as \posit<n,2>. Often the ``2'' is dropped; for example, the standard 64-bit posit, \posit<64,2>,  would be written as \p<64>.

For example, Figure \ref{fig:posit} shows the bitstring of the number ${3.5465 \times 10^{-6}}$ encoded as a \posit<16,2>.
In this example, we have a run length of $r = 5$, therefore $k = -5$.  
The exponent  $e = 1$ and $m = 7$ fraction bits.


\begin{figure}[ht!]

\begin{center}
\begin{tikzpicture}[scale=1.0]

\definecolor{sign}{RGB}{97, 160, 225}
\definecolor{fraction}{RGB}{255, 196, 46}  
\definecolor{expo}{RGB}{134,200,79}
\definecolor{regime}{RGB}{255,122,35}



\def\fracbits{7}   
\def\ebits{2}       
\def\regbits{6}

\def\regstr{0/1, 0/2, 0/3, 0/4, 0/5, 1/6}
\def\signstr{0/1}                       
\def\expostr{0/1, 1/2}   
\def\fracstr{1/1, 1/2, 0/3, 1/4, 1/5, 1/6, 0/7}  

\def\offset{0.7}   
\def\step{0.5}      
\def\y{2}           

\pgfmathsetmacro{\nbits}{\fracbits + \ebits + \regbits + 1}
\pgfmathsetmacro{\exponent}{\ebits * \step}
\pgfmathsetmacro{\regime}{\regbits * \step}
\pgfmathsetmacro{\mantissa}{\fracbits * \step}
\pgfmathsetmacro{\x}{\step + \offset}
\pgfmathsetmacro{\yy}{\step + \y}

\draw[step=\step cm, black, very thin] (0,\y-0.01) grid (\nbits*\step,\yy);

\fill[sign, fill opacity=.4] (0,\y-0.01) rectangle (\step, \yy);

\draw[] (0,\y+\offset) -- ++ (0,0.50em) -| (\step,\y+\offset);
\node (s) at (\step / 2,\y+1.7*\offset) {\small\textsf{sign}};

\fill[regime, fill opacity=.4] (\step,\y-0.01) rectangle (\step + \regime, \yy);  

\draw[] (\step,\y+\offset) -- ++ (0,0.50em) -| (\step + \regime,\y+\offset);
\node (r) at (\regime / 2 + \step,\y+1.7*\offset) {\small\textsf{regime}};

\fill[expo, fill opacity=.4] (\step + \regime,\y-0.01) rectangle (\step + \regime + \exponent, \yy);  

\draw[] (\step + \regime,\y+\offset) -- ++ (0,0.50em) -| (\step + \regime + \exponent,\y+\offset);

\node (e) at (\regime + \exponent / 2 + \step,\y+1.7*\offset) {\small\textsf{exp}};

\fill[fraction, fill opacity=.4] (\step + \regime + \exponent,\y-0.01) rectangle ( \step + \regime + \exponent + \mantissa, \yy);

\draw[] (\exponent+\step + \regime, \y+\offset) -- ++ (0,0.50em) -| (\step + \regime + \exponent + \mantissa,  \y+\offset);
\node (f) at (\step + \regime + \mantissa/2 + \exponent,\y+ 1.7*\offset) {\small\textsf{fraction (\fracbits \, bits)}};

\filldraw[black] (4.5,2.25) circle (1.75pt);


\foreach\i\j in \signstr
  {
  \node (s\i) at (\j*\step-\step/2, \yy-\step/2) {\texttt{\i}};
  } 

\foreach\i\j in \regstr
  {
  \node (r\i) at (\j*\step + \step/2, \yy-\step/2) {\texttt{\i}};
  }

\foreach\i\j in \expostr
  {
  \node (e\i) at (\regime + \j*\step + \step/2, \yy-\step/2) {\texttt{\i}};
  }

\foreach\i\j in \fracstr
  {
  \node (f\i) at (\exponent + \regime + \j*\step + \step/2, \yy-\step/2) {\texttt{\i}};
  }

\end{tikzpicture}
\end{center}

 \caption{A \posit<16,2> with 7-bit fraction and 2-bit exponent.}
    \label{fig:posit}
\end{figure}
The decoding given by (\ref{eqn:decode}) is: 
\begin{equation}
   \left( \,(1 - 3{\cdot 0})  +{\frac{1}{2}} + {\frac{1}{4}} + {\frac{0}{8} } + {\frac{1}{16} } + {\frac{1}{32} } + { \frac{1}{64}} + { \frac{0}{128}} \right) \times 2^{ 4(-5) + 1} .
\end{equation}
A posit visualization tool is available online \cite{positviz}.

The posit standard specifies using an accumulation register called a \textit{quire} to defer rounding. 
According to the Standard \cite[p. 6]{posits2022}, ``a quire value is either \texttt{NaN} or an integer multiple of the square of the minimum positive value, represented as a 2's complement binary number with $16n$ bits."  
There are 256 bits reserved for this accumulator in a \posit<16,2> configuration.
A single round is performed to produce the final representation.

Table \ref{tab:positprops} lists properties, including the maximum consecutive integer, the range of fraction lengths, quire precision, and the minimum and maximum positive representable number in the posit standard.
There are always $2^{n-2}+1$ values in the interval $[0,1]$, regardless of the number of exponents bits.
Table \ref{tab:comparisons} displays comparisons between various posit configurations and floating-point formats.

%

\begin{table}
\centering
     \caption{Posit Properties under the standard $E = 2$.} \label{tab:positprops}
    \begin{tabular}{l l l}
    \hline
    Property & \;\;\;\;\; & Value \\
    \hline 
        Fraction length & &  $0$ to $\max(0,n-5)$ bits\\ 
        Minimum Positive Value  & & $2^{-4n+8}$ \\
         Maximum Positive Value  & & $2^{4n - 8}$ \\
         Maximum Consecutive Integer   &   & $\lceil{2^{\lfloor{4(n-3)/5}\rfloor} }\rceil$ \\         
         Quire format precision & & $16n$ bits \\
         Quire sum limit & & $2^{23 + 4n}$ \\
         Number of values in $[0,1]$ & & $2^{n-2}+1$\\
         \hline
    \end{tabular}
\end{table}

\begin{table}[htbp]
    \setlength{\tabcolsep}{6pt}
 \caption{Specifications of multiple binary formats for comparison including IEEE (half, single, and double precision), Google's {\tt bfloat16}, and three standard posit configurations.  The precision of the arithmetic is measured by the unit round-off, listed below as the distance from 1.0 to the next larger representable value.}   \centering
    \begin{tabular}{lccccc}
            \hline
    \textbf{Format}          & \textbf{Fraction}     & \textbf{Exponent}     & \textbf{u}     & \textbf{min}     & \textbf{max}   \\
            \hline
    {posit$\langle$16,2$\rangle$}     & $11^*$ & 2 & $\num{2.44e-04}$  & $\num{1.39e-17}$  & $\num{7.21e+16}$   \\
    {posit$\langle$32,2$\rangle$}     & $27^*$ & 2 & $\num{7.45e-09}$  & $\num{7.52e-37}$  & $\num{1.33e+36}$   \\
    {posit$\langle$64,2$\rangle$}     & $59^*$ & 2 & $\num{1.73e-18}$  & $\num{2.21e-75}$  & $\num{4.52e+74}$   \\
    {bfloat16} & $8$ & $8$ & $\num{3.91e-03}$ & $\num{1.18e-38}$ & $\num{3.39e+38}$ \\ 
    {fp16}         & 11  & 5  & $\num{4.88e-04}$   & $\num{6.10e-05}$  & $\num{6.55e+04}$   \\
     {fp32}         & 24  & 8  & $\num{1.19e-07}$   & $\num{1.18e-38}$  & $\num{3.40e+38}$   \\
    {fp64}         & 53  & 11  & $\num{1.11e-16}$   & $\num{2.22e-308}$  & $\num{1.80e+308}$   \\
        \hline
        \multicolumn{6}{l}{$^*$ \small denotes variable fraction bits (maximum number of bits listed)}
    \end{tabular}
  \label{tab:comparisons}
\end{table}

%
\section{Numerical Experiments}\label{sec:numx}
We perform numerical experiments to compare the convergence of the solution to the system $Ax = b$ using mixed-precision iterative refinement (Algorithm \ref{alg:mixir}) with the posit number system.  
Three algorithms converting to low-precision posits will be examined.
For all systems, the exact solution is $x = (1, 1, \dots, 1)^T$.

\subsection{Posit Configurations. } The incorporation of the quire in our experiments facilitates the utilization of two posit configurations: \posit<32,2> for maintaining working precision and \posit<16,2> for conducting low-precision calculations. 
This approach contrasts with the requisite use of three precisions in analogous experiments performed within the IEEE floating-point framework \cite{haidar2018design,haidar2018harnessing,haidar2017investigating,higham2019new,higham2019squeezing}.

\subsection{Test-Suite. } 
Table \ref{tab:mats} lists ten nonsingular square matrices extracted from the \textit{SuiteSparse Matrix Collection}\footnote{Formerly known as the University of Florida Sparse Matrix Collection.} \cite{davis2011university} sorted by the unique identifier (ID), along with key characteristics such as the size, percent of nonzero entries, maximum and minimum absolute (nonzero) entry, and condition number.
We restrict the matrices' size to $n \le 240$ due to the time required for LU factorization without native hardware support for posits.
These matrices represent problems in diverse fields such as computational fluid dynamics (CFD), chemical simulation, materials science, optimal control, structural mechanics, and 2D/3D sequencing. 
The matrices were chosen to facilitate a (near) comparative analysis with Higham et al.\,\cite{higham2019squeezing}, (and others), who employed \fp<16> for the low-precision computations,  \fp<32> as the working precision, and \fp<64> for residual calculations in the iterative refinement process.
However, one difference is nearly all the entries lie within the dynamic range of \low<16>.

\begin{table}[ht!]
\setlength{\tabcolsep}{9pt}
 \caption{Selected matrices representing various disciplines from the Suite Sparse Matrix Collection.  
 ID is the matrix identifier as listed in the collection.  A `$^*$' denotes a symmetric matrix.}
  \centering
    \begin{tabular}{llccccc}
    \hline
    \textbf{ID} & \textbf{Matrix} & $n$  & nnz(\%) & $\max {|a_{ij}|}$  & $\min {|a_{ij}| \ne 0}$ &  $\kappa_\infty(A)$ \\
    \hline
    6 & arc130 & 130 & 6.14& $\num{1.05e+05}$ & $\num{7.71e-31}$ & $\num{1.20e+12}$\\ 
    23$^*$ & bcsstk01 & 48 & 17.36& $\num{2.47e+09}$ & $\num{3.33e+03}$ & $\num{1.60e+06}$\\ 
    27 & bcsstk05 & 153 & 10.35 & $\num{1.43e+04}$ & $\num{4.65e-10}$ & $\num{3.53e+04}$\\ 
    206$^*$ & lund\_a & 147 & 11.33 & $\num{1.50e+08}$ & $\num{1.22e-04}$ & $\num{5.44e+06}$\\ 
    217$^*$ & nos1 & 237 & 1.81  & $\num{1.22e+09}$ & $\num{8.00e+05}$ & $\num{2.53e+07}$\\ 
     232 & pores\_1 & 30 & 20.00 & $\num{2.46e+07}$ & $\num{4.00e+00}$ & $\num{2.49e+06}$\\ 
     239 & saylr1 & 238 & 1.99 & $\num{3.06e+08}$ & $\num{7.19e-04}$ & $\num{1.59e+09}$\\ 
    251 & steam1 & 240 & 3.90  & $\num{2.17e+07}$  & $\num{1.48e-07}$ & $\num{3.11e+07}$  \\ 
    263 & west0132 & 132 & 2.37   & $\num{3.16e+05}$  & $\num{3.31e-05}$ & $\num{1.05e+12}$  \\ 
    298 & bwm200 & 200 & 1.99  & $\num{6.15e+02}$  & $\num{4.00e+00}$ & $\num{2.93e+03}$  \\ 
    \hline
    \end{tabular}
  \label{tab:mats}
\end{table}

\subsection{Converting to Low Precision Posits.}  
Many matrices have entries outside the dynamic range of these low-precision configurations.
In IEEE half (or single) precision, the elements of $A$ may overflow or underflow during the transition to lower precision.
However, with posits, matrix entries exceeding the dynamic range will be mapped to the maximum or minimum signed value, thus preventing overflow or underflow.
Algorithm \ref{alg:21} converts a higher precision posit to a lower precision.

%
\begin{algorithm}[ht]
\caption{Round to low-precision posit}\label{alg:21}
\begin{algorithmic}[1] 
\Require $A \in \mathbb{R}^{n \times n}$
\Ensure $A^{(L)}$ lower configuration 

\State $A^{(L)} = \rm{fl}(A)$ using \posit<16,2> 

\State $\forall i, j$ such that $|a_{ij}^{(L)}| > x_{\rm{max}}$, set $a_{ij}^{(L)} = \text{sgn}(a_{ij}) \, x_{\rm{max}}$
\State $\forall i, j$ such that $|a_{ij}^{(L)}| < x_{\rm{min}}$, set $a_{ij}^{(L)} = \text{sgn}(a_{ij}) \, x_{\rm{min}}$
\end{algorithmic}
\end{algorithm}

The second approach, Algorithm \ref{alg:22}, follows the method proposed in  \cite{higham2019squeezing}. 
In this approach, the matrix is scaled before rounding.
After scaling and rounding the matrix, Algorithm \ref{alg:mixir} will be applied.
The right-hand side vector $b$ must also be scaled by the same factor $\mu$ to preserve the equality of the system after scaling the matrix.

\begin{algorithm}[H]
\caption{Scale matrix entries, then round to low-precision}\label{alg:22}
\begin{algorithmic}[1]  
  \Require $A \in \mathbb{R}^{n \times n}$
  \Ensure $A^{(L)}$ low-precision    
  \State $A = \mu A$ \;\;\; for $\mu \in \mathbb{R}$
  \State $A^{(L)} = {\rm{fl}}(A)$
\end{algorithmic}
\end{algorithm}

The third conversion strategy uses two-sided diagonal scaling.
Since the condition number $\kappa(A) = ||A|| \cdot ||A^{-1}||$ affects error bounds, for $\kappa(A) \gg 1$, iterative refinement may never converge \cite{carson2017new}. 
To mitigate this issue, the system is often rescaled to reduce the condition number, aiming to improve convergence. 
However, it is crucial to note that scaling strategies can be unreliable, as they may produce more ill-conditioned matrices \cite{forsythe1967computer}. 
Therefore, it is generally recommended to apply scaling strategies on a problem-by-problem basis, taking into consideration the significance of the entries $a_{ij}$, the units of measurement, and error in the data \cite{golub2013matrix}.

This work's primary technique is applying simple two-sided row and column scalings to the matrix $A$. 
These scalings aim to equilibrate $A$, as described in the references \cite{forsythe1967computer} and \cite{golub2013matrix} to reduce the condition number of $A$.

In particular, equilibration works to solve $Ax = b$ by first mutating the system by premultiplying by the diagonal matrix $D^{-1}_1$ to get  
$
D^{-1}_1A x = D^{-1}_1b.
$
The goal is to choose $D^{-1}_1$ so that $D_1A$ has approximately the same infinity norm, thus making it less likely to combine numbers of extreme opposite scale during the elimination process \cite{golub2013matrix}.
Next, 
let $x = D_2y$, and solve the modified system, 
\[
D^{-1}_1 A D_2 y = D^{-1}_1 b, 
\]
where $D_2$ is a diagonal matrix that scales the columns.  
\[
D^{-1}_1 A D_2 y = D^{-1}_1 b \Rightarrow y  = D_2^{-1}A^{-1} b. 
\]
After post-processing, we have
\[
D_2 y  = A^{-1} b = x .
\]
Algorithm \ref{alg:23} outlines a basic implementation of this two-sided diagonal scaling process.  
LAPACK \cite{anderson1999lapack} uses the nearest power of 2 scaling to avoid rounding errors (see \texttt{xyyEQUB} routines).
It is important to note that the order in which scaling is applied matters (e.g., column first scaling), resulting in different matrices with different condition numbers (see examples and discussion \cite[p.45]{forsythe1967computer}).
 
\begin{algorithm}[H]
\caption{(Row and Column Equilibration).  Compute diagonal matrices $D_1$ and $D_2$ such that $D^{-1}_1 A D_2$ has a maximum element of 1 in any row and column.}\label{alg:23}
\begin{algorithmic}[1]
\Require $A \in \mathbb{R}^{n \times n}$
\Ensure $A^{(L)}$ 
\State $D_1 = I$
\For{$i = 1:n$}
    \State {$D_1(i,i) = \|A(i,:)\|_\infty$}
\EndFor
\State $B = D_1A$
\State $D_2 = I$
\For{$j = 1:n$}
    \State {$D_2(j,j) = \|B(:,j)\|_\infty^{-1}$}
\EndFor
\State $A^{(L)} = {\rm{fl}}(\mu D^{-1}_1AD_2)$ 
\end{algorithmic}
\end{algorithm}

\subsection{Algorithm performance. }Algorithm performance is measured by the total number of iterations required to converge.
Convergence is assessed using \textit{Criterion 1} in  \cite[p.54]{barrett1994templates}, where the value of $\epsilon = 10^{-8}$ is set, such that
\begin{align}
    \| r_i \| \le \epsilon \left(\|A\| \cdot \| x_i \| + \| b \| \right)
\end{align}
yielding the forward error bound
\[
\| x - x_i \| \le \epsilon \|A^{-1}\|  \left(\|A\| \cdot \| x_i \| + \| b \| \right).
\]
To determine the effectiveness of scaling on reducing the condition number, we approximate $\|A^{-1}\| \approx \|y\| / \|x\|$, where $Ay = x$, following \cite{cline1979estimate}.

\subsection{Software.} 
We utilized the \textit{Universal Numbers Library} \cite{omtzigt2023}, or \textit{Universal} for short. 
\textit{Universal} is a C++ header-only template library that implements multiple machine number representations, including floats and posits, with arbitrary configurations and supports standard arithmetic operations. Clang 14 (clang-1400.0.29.202) was used to compile the code.

\subsection{Hardware.} 
All numerical experiments were conducted on a desktop computer with a 3.2 GHz 6-Core Intel Core i7 processor, 32 GB of 2267 MHz DDR4 RAM, and a macOS 13.5 operating system.


\section{Results}\label{sec:results}
The results are displayed in Table \ref{tab:results}. 
The method achieved convergence across all test cases by employing Algorithm \ref{alg:23} to cast elements to a lower precision.
However, neither Algorithms \ref{alg:21} and \ref{alg:22} guarantee convergence.  
In general, two-sided scaling results in the best performance using iterative refinement.
We highlight the approaches that involved scaling and conversion, operating without preconditioning strategies like incomplete LU factorization or preconditioned GMRES.
We found convergence with various values of $\mu$, for example, when $\mu = 16$, Algorithm \ref{alg:23} converged for nos1, and with $\mu = {0.75}$, Algorithm \ref{alg:22} converged for saylr1 and west0132, and bwm200.
Algorithm \ref{alg:21} converged for all test matrices using  \posit<24,2>.

\begin{table}[ht!]
\setlength{\tabcolsep}{12pt}
 \caption{Total number of iterations required by the mixed-precision iterative refinement method with scaling factor $\mu = 1/16$. 
 Divergence is indicated by $\infty$. 
 Convergence was obtained for other values of $\mu$ in Algorithm \ref{alg:22}. 
 For the case ``nos1'', indicated by `$^*$', convergence was achieved with $\mu=16$.}\label{tab:results}
  \centering
    \begin{tabular}{lccc}
    \hline
    \textbf{Matrix}  & \textbf{Algorithm \ref*{alg:21}} & \textbf{Algorithm \ref*{alg:22}} & \textbf{Algorithm \ref*{alg:23}} \\
    \hline
    arc130  &  2   &   1   &   1\\
    bcsstk01  &  $\infty$   &   $6$    &   4\\
    bcsstk05  &  $\infty$   &   $5$   &   8\\
    lund\_a  &  $\infty$   &   $27$   &   5\\
    nos1  &  $\infty$   &   $180$   &   $81^*$\\
    pores\_1  &  $14$   &   5   &   3\\
    saylr1  &  $\infty$   &   $\infty$   &   87\\
    steam1  &  $\infty$   &   $\infty$   &   2\\
    west0132  &  $\infty$   &   $\infty$   &   4\\
    bwm200  &  $\infty$   &   $63$   &   9\\
    \hline 
\end{tabular}
\end{table}

Table \ref{tab:scaledcon} illustrates the reduction in condition numbers resulting from the implementation of Algorithms \ref{alg:22} and \ref{alg:23}, which further explores the impact of scaling. 
The condition numbers displayed are computed with working precision. 
The table indicates that diagonal scaling reduces large condition numbers for the evaluated matrices. 
Notably, for matrices 27 and 298, the scaled matrices remain ill-conditioned, emphasizing that convergence is not solely determined by the condition number, as mentioned in \cite{golub2013matrix}.
Additionally, Table \ref{tab:unique} presents the count of distinct matrix element values both before and after the application of Algorithms \ref{alg:21} through \ref{alg:23}. 
The table emphasizes how these algorithms' scaling procedures affect the matrix elements' variability.

\begin{table}[ht!]
\setlength{\tabcolsep}{12pt}
 \caption{Condition number after application of Algorithm.}\label{tab:scaledcon}
  \centering
   \begin{tabular}{lcccc}
   \hline
   &   & Algorithm \ref{alg:22} & Algorithm \ref{alg:23} \\
   \cline{3-4}
Matrix  & $\kappa_\infty(A)$ & $\kappa_\infty(B)$ & $\kappa_\infty(B)$ \\
\hline
arc130   & $\num{1.20e+12}$     & $\num{1.20e+12}$   & $\num{4.07e+02}$   \\
bcsstk01  & $\num{1.60e+06}$   &  $\num{1.59e+06}$   & $\num{3.25e+03}$   \\
bcsstk05  & $\num{3.53e+04}$     & $\num{1.43e+04}$   & $\num{1.43e+04}$   \\
lund\_a& $\num{5.44e+06}$   &  $\num{5.44e+06}$   & $\num{8.29e+04}$   \\
nos1 & $\num{2.53e+07}$     & $\num{2.53e+07}$   & $\num{9.07e+06}$   \\
pores\_1 & $\num{2.49e+06}$  & $\num{2.49e+06}$   & $\num{8.20e+03}$   \\
saylr1 & $\num{1.59e+09}$   &  $\num{1.59e+09}$   & $\num{1.96e+05}$   \\
steam1  & $\num{3.11e+07}$   & $\num{3.11e+07}$   & $\num{9.27e+00}$   \\
west0132 & $\num{1.05e+12}$   &  $\num{1.05e+12}$   & $\num{9.59e+06}$   \\
bwm200 & $\num{3.03e+03}$   & $\num{3.03e+03}$    & $\num{2.37e+03}$  \\
\hline 
\end{tabular}
\end{table}

\begin{table}[ht!]
\setlength{\tabcolsep}{12pt}
 \caption{Number of unique matrix elements after applying Algorithms \ref{alg:21} -- \ref{alg:23}.  The number of unique entries of the original matrix before conversion is denoted in parentheses.}\label{tab:unique}
  \centering
    \begin{tabular}{lccc}
    \hline
    \textbf{Matrix}  & \textbf{Algorithm \ref*{alg:21}} & \textbf{Algorithm \ref*{alg:22}} & \textbf{Algorithm \ref*{alg:23}} \\
    \hline
    arc130 (961)   &    392  &  224     &  339\\
    bcsstk01 (105) &    15   &  103     &  232\\
    bcsstk05 (246) &    121  &  244     &  725\\
    lund\_a (232)  &     74  &  28      &  74\\
    nos1  (9) &  3   &   9   &   7\\
    pores\_1 (152) &  102   &   131   &   132\\
    saylr1 (705) &  453   &   277   &   665\\
    steam1 (764) &  592   &   325   &   767\\
    west0132 (237) &  227   &   145   &   167\\
    bwm200 (7) &  7   &   7   &   6\\
    \hline 
\end{tabular}
\end{table}

\section{Summary and Conclusions}\label{sec:conclusions}

We utilized mixed-precision iterative refinement with posits to accurately approximate the solution to the equation $Ax = b$ for ten different sparse matrices from the SuiteSparse Matrix Collection. We employed a deferred rounding mechanism in the software to simulate high-precision residual calculations. After applying three conversion algorithms, we observed that employing two-side diagonal scaling resulted in convergence in all cases and a significant reduction in the condition number in most cases without using a dedicated preconditioner.
 Our experiments showed convergence similar to those using \fp<16>. For most matrices in the test suite, one-sided row scaling produced similar results. The convergence and convergence rate generally depends on the matrix type, condition number, and matrix size \cite{haidar2017investigating}. However, we found that the matrix size was not a determining factor in our experiments. Our findings extended the work of \cite{higham2019squeezing} to posits and \cite{buoncristiani2020evaluating} to nonsymmetric matrices using LU factorization instead of Cholesky.

Converting matrix entries to a lower precision for solving $Ax = b$ is a crucial yet challenging task in numerous applications, such as training deep learning models \cite{carmichael2019deep,lu2020evaluations,murillo2020deep}, image processing \cite{cococcioni2020fast}, and climate modeling \cite{dawson2018reliable}.
In addition to the loss of precision during conversion to low precision, which typically produces overflow or underflow, values are mapped to the maximum and minimum values in the range by default, which can produce undesirable behaviors.
While \fp<16> has become a popular choice due to the increased availability of hardware support, there is a need to explore alternative formats, such as posits, in which hardware is becoming available \cite{mallasen2022percival}.
Scaling strategies that leverage the distribution of machine numbers in the posit number system need to be further understood.

This research establishes a foundation for exploring various directions and conducting further analyses.
With posits, techniques that couple iterative refinement with preconditioned GMRES and incorporate incomplete or block low-rank LU factorization, among others, need to be investigated. 
While optimal scaling techniques have been discussed in \cite{bauer1963optimally}, determining when and how to scale effectively remains an open problem \cite{elble2012scaling}. Various scaling strategies exist, including Hungarian scaling \cite{hook2019max} and simple interval transformation \cite{higham2019squeezing}, which could be explored in the context of posits. 
Lastly, new number systems such as \textit{takums} developed by \cite{hunhold2024beating} offer further directions for research related to low and mixed-precision algorithms.

\bibliographystyle{ams}
\bibliography{references}

\end{document}